# Symbolic and High-Accuracy Solutions to Differential and Integral Problems via a Novel Recursive Inverse Laplace Method

Mohamed Mostafa
Physics and engineering mathematics department, Faculty of Engineering - Helwan University – Mataria Branch, Cairo, Egypt.
Email: mohamed_mostafa@m-eng.helwan.edu.eg

**Abstract**

In this paper, we introduce a novel semi-analytical method for solving a broad class of initial value problems involving differential, integro-differential, and delay equations, including those with fractional and variable-order derivatives. The proposed approach is based on the inverse Laplace transform, applied initially—unlike traditional Laplace-based techniques which begin with a forward transformation. By assuming the unknown solution is the Laplace transform of an auxiliary function, the method reformulates the problem in the time domain and recursively solves for this function using symbolic operations. The final solution is then obtained by applying the Laplace transform to the result. This strategy enables the construction of symbolic solutions as generalized logarithmic-power series with arbitrary accuracy, and naturally accommodates complex terms such as $x^{\alpha} \ln(x)^{n}$. The method is highly versatile and demonstrates superior speed and precision across a wide range of linear and nonlinear problems, including singular, fractional, and chaotic systems. Several benchmark examples are provided to validate the reliability and efficiency of the proposed technique compared to classical numerical methods. The results confirm that the new method offers a powerful and flexible framework for symbolic computation of initial value problems.

**Keywords** Recursive inverse Laplace transform, Symbolic solution, Initial value problems, Delay and integro-differential equations, Arbitrary accuracy, Variable-order systems.

## 1. Introduction

Initial value problems (IVPs) involving differential, integro-differential, and delay equations are foundational in modeling phenomena in engineering, physics, biology, and control systems. These models often include intricate features such as memory effects, singularities, time delays, or fractional and variable-order derivatives. Solving such problems analytically or symbolically remains a central challenge in applied mathematics.

Analytical and semi-analytical techniques like the Adomian Decomposition Method (ADM), Homotopy Analysis Method (HAM), and Variational Iteration Method (VIM) have been widely used. They are typically combined with the Laplace transform to convert differential equations into algebraic forms, leading to hybrid schemes such as the Laplace–Adomian Decomposition Method (LADM) and Laplace–Homotopy Perturbation Method (LHPM). Zhang and Wang [1] applied LADM to fractional Fokker–Planck and Boussinesq–Burger equations, showing its efficiency and simplicity—and achieving excellent agreement with exact solutions. Elbadri [2] introduced a generalized LADM for IVPs with generalized fractional derivatives, demonstrating rapid convergence in symbolic series. Shah et al. [3] applied LADM to third-order dispersive fractional PDEs and achieved highly accurate convergent series solutions. Eltayeb et al. [4] employed multi-dimensional LADM for singular fractional pseudo-hyperbolic equations, yielding accurate symbolic and approximate results in 2D examples.

However, classical Laplace-based approaches encounter notable limitations. These include difficulty handling solutions involving non-polynomial terms such as $x\ln(x)$, inefficiencies when addressing variable-order derivatives, and relying on complex structures like Adomian polynomials, which increase computational cost. These limitations have motivated efforts to develop symbolic and more structurally transparent solution techniques.

A major recent development has been made by Dana Mazraeh and Parand [5] who designed a method that fuses deep Q-learning with symbolic grammar-based exploration, generating highly structured symbolic solutions for nonlinear differential systems.

Despite these advances, most symbolic or hybrid methods either depend on external training processes or generate solutions through indirect symbolic matching. In contrast, we propose a new semi-analytical method—the Recursive Inverse Laplace Transform (RILT)—which reverses the classical Laplace framework by applying the inverse transform first. The central idea is to assume that the unknown solution is the Laplace transform of an auxiliary function. This enables the original equation to be converted into a recursive symbolic form directly in the time domain.

This inverse-first strategy facilitates the symbolic construction of generalized logarithmic-power series, allowing the method to tackle singularities, variable-order operators, and nonlinearities naturally. The method supports symbolic solutions of arbitrary accuracy and parameter dependence without relying on Adomian polynomials or training-based AI.

The remainder of the paper is organized as follows. Section 2 introduces the necessary preliminaries and definitions related to fractional calculus and Laplace transforms. Section 3 presents key lemmas that form the theoretical foundation of the proposed method. In Section 4, we describe the algorithm in detail and discuss its implementation. Section 5 demonstrates the effectiveness of the method through a diverse set of benchmark problems. Finally, conclusions and future research directions are discussed in Section 6.

## 2. Preliminaries

This section provides essential definitions and properties required to develop the proposed method. Specifically, we recall the Dirac delta function, the Caputo fractional derivative of variable order, and the Laplace transform. These concepts are foundational for the symbolic and recursive framework adopted in this work.

### 2.1. Dirac Delta Function [6]

The Dirac delta function, denoted by $\delta(t)$, is a distribution defined such that

$$\delta(t) = \begin{cases} 0, & t \neq 0, \\ \infty, & t = 0, \end{cases} \quad \text{with} \quad \int_{-\infty}^{\infty} \delta(t)dt = 1.$$

It satisfies the sifting property:

$$\int_a^b f(t)\delta(t)dt = f(0), \qquad 0 \in [a, b],$$

and for its $n$-th distributional derivative $\delta^{(n)}(t)$, we have:

$$\int_a^b f(t)\delta^{(n)}(t)dt = (-1)^n f^{(n)}(0), \qquad 0 \in [a, b]. \tag{1}$$

The Laplace transform of a shifted derivative of the delta function is given by:

$$\mathcal{L}\{\delta^{(n)}(t-a)\} = s^n e^{-as}, \qquad a \geq 0. \tag{2}$$

### 2.2. Caputo Fractional Derivative of Variable Order [6]

Let $\alpha(x)$ be a variable order function with $n - 1 < \alpha(x) < n$, where $n \in \mathbb{Z}^+$. The Caputo fractional derivative of a function $f(x)$ of order $\alpha(x)$ is defined as:

$$D^{\alpha(x)}f(x) = \begin{cases} \frac{1}{\Gamma(n-\alpha(x))}\int_0^x \frac{f^{(n)}(\tau)}{(x-\tau)^{\alpha(x)-n+1}}d\tau, & \alpha(x) \notin \mathbb{N}, \\ f^{(m)}(x), & \alpha(x) = m \in \mathbb{N}. \end{cases}$$

For the power function $f(x) = x^p$, the Caputo fractional derivative of variable order is given by:

$$D^{\alpha(x)}x^p = \Theta(p-\alpha(x)) \cdot \frac{\Gamma(p+1)}{\Gamma(p-\alpha(x)+1)}x^{p-\alpha(x)}, \qquad p \in \mathbb{R}^+, \qquad \alpha(x) \geq 0,$$

where $\Theta(\cdot)$ is the Heaviside step function.

**2.3. Laplace Transform** [4]

The Laplace transform of a function $f(t)$, defined on $[0,\infty)$, is given by:

$$\mathcal{L}\{f(t)\} = F(s) = \int_0^\infty f(t)e^{-st}dt, \qquad s > 0.$$

The transformation exists for piecewise continuous functions of exponential order $\beta$, i.e. $|f(t)| \leq Me^{\beta t}$ for some constants $\beta, M > 0$.

Key properties of the Laplace transform include:

- Linearity:

$$\mathcal{L}\{a_1 f(t) + a_2 g(t)\} = a_1\mathcal{L}\{f(t)\} + a_2\mathcal{L}\{g(t)\}, \qquad a_1, a_2 \in \mathbb{R}.$$

- Translation in time:

$$\mathcal{L}\{\Theta(t-a)f(t-a)\} = e^{-as}\mathcal{L}\{f(t)\}, \qquad a \geq 0.$$

- The inverse Laplace transform is formally defined via the Bromwich integral:

$$f(t) = \mathcal{L}^{-1}\{F(s)\} = \frac{1}{2\pi i}\lim_{T\to\infty}\int_{\gamma-iT}^{\gamma+iT} e^{st}F(s)ds,$$

where $\gamma$ is a real constant such that the contour of integration lies within the region of convergence of $F(s)$.

## 3. Fundamental Lemmas

This section presents three lemmas that form the theoretical backbone of the proposed method. They provide the link between symbolic expressions in the Laplace domain and their recursive representations in the time domain, especially in the presence of logarithmic-power terms and variable-order derivatives.

**Lemma 1.**

Let $u(x) = \int_0^\infty f(t)x^t dt$, where $x = e^{-s}$, and assume $f(t)$ is sufficiently smooth. Then the Caputo fractional derivative of $u(x)$ of order $\alpha(x)$ satisfies:

$$x^{\alpha(x)} D^{\alpha(x)} u(x) = \mathcal{L}\left\{\Theta\big(t - \alpha(x)\big) \cdot \frac{\Gamma(t+1)}{\Gamma(t - \alpha(x) + 1)} f(t)\right\}.$$

**Proof.**

By assuming $u(x) = \mathcal{L}\{f(t)\} = \int_0^\infty f(t)x^t dt$, and applying the Caputo derivative under the integral sign using the fractional power rule,

$$D^{\alpha(x)} u(x) = \int_0^\infty f(t) D^{\alpha(x)} x^t dt,$$

and from the Caputo derivative of $x^t$, we get:

$$D^{\alpha(x)} x^t = \Theta\big(t - \alpha(x)\big) \cdot \frac{\Gamma(t+1)}{\Gamma(t - \alpha(x) + 1)} x^{t - \alpha(x)}.$$

Multiplying both sides by $x^{\alpha(x)}$ yields the desired result:

$$\begin{aligned} x^{\alpha(x)} D^{\alpha(x)} u(x) &= \int_0^\infty \left(\Theta\big(t - \alpha(x)\big) \cdot \frac{\Gamma(t+1)}{\Gamma(t - \alpha(x) + 1)} f(t)\right) x^t dt \\ &= \mathcal{L}\left\{\Theta\big(t - \alpha(x)\big) \cdot \frac{\Gamma(t+1)}{\Gamma(t - \alpha(x) + 1)} f(t)\right\}. \blacksquare \end{aligned}$$

**Lemma 2.**

Let $f(t)$ be differentiable and $x = e^{-s}$. Then,

$$\mathcal{L}\big\{f(t) \cdot \mathcal{L}^{-1}\{x^k \ln(x)^m\}\big\} = \frac{\partial^m}{\partial t^m}\big(f(t)x^t\big)\Big|_{t=k}, \qquad m \in \mathbb{N}, \qquad k \geq 0.$$

**Proof**.

By using (1) and (2) where $x = e^{-s}$,

$$\mathcal{L}\big\{f(t) \cdot \mathcal{L}^{-1}\{x^k \ln(x)^m\}\big\} = \int_0^\infty f(t) \cdot (-1)^m \delta^{(m)}(t-k) x^t dt = \frac{\partial^m}{\partial t^m}\big(f(t)x^t\big)\Big|_{t=k}. \blacksquare$$

**Lemma 3.**

$$D^{\alpha(x)}\big(x^k \ln(x)^m\big) = \frac{\partial^m}{\partial t^m}\left(\Theta\big(t - \alpha(x)\big) \frac{\Gamma(t+1)}{\Gamma\big(t + 1 - \alpha(x)\big)} x^{t - \alpha(x)}\right)\Bigg|_{t=k}, m \in \mathbb{N}, \qquad k \geq 0.$$

**Proof**.

By using (1) and (2) where $x = e^{-s}$,

$$
\begin{aligned}
D^{\alpha(x)}(x^k \ln(x)^m) &= D^{\alpha(x)} \int_0^\infty (-1)^m \delta^{(m)}(t-k) x^t dt \\
&= \int_0^\infty (-1)^m \delta^{(m)}(t-k) \Theta\big(t-\alpha(x)\big) \frac{\Gamma(t+1)}{\Gamma\big(t+1-\alpha(x)\big)} x^{t-\alpha(x)} dt \\
&= \frac{\partial^m}{\partial t^m} \left( \Theta\big(t-\alpha(x)\big) \frac{\Gamma(t+1)}{\Gamma\big(t+1-\alpha(x)\big)} x^{t-\alpha(x)} \right)\Bigg|_{t=k}. \blacksquare
\end{aligned}
$$

# 4. Proposed Method

This section presents the development of the proposed method for solving IVPs involving differential, integro-differential, and delay equations, including those with variable or fractional order derivatives. The method is based on assuming that the unknown solution is the Laplace transform of an auxiliary function. By reversing the classical Laplace-based solution strategy, the proposed method enables symbolic, recursive construction of solutions with arbitrary accuracy.

## 4.1. Method Overview

Let the general form of an IVP be written as:

$$x^{\alpha(x)} D^{\alpha(x)}\big(u(x)\big) = A(u(x), x), \tag{3}$$

with the initial conditions

$$u^{(k)}(0) = c_k, \qquad k = 0,1,\dots,n-1, \qquad n-1 < \alpha(x) < n,$$

where $A(\cdot)$ is the integral/differential operator applied to an expression of $u(x)$, and $\alpha(x)$ is the highest order derivative in the equation.

Unlike conventional approaches that apply the Laplace transform directly to the equation, we assume the solution $u(x)$ is itself a Laplace transform:

$$u(x) = \mathcal{L}\{f(t)\} = \int_0^\infty f(t) x^t dt.$$

The key steps of the method are:

1. Apply the inverse Laplace transform to both sides of the original equation.
2. Use lemmas established in Section 3 to express the equation in terms of $f(t)$ in the time domain.
3. Solve recursively or symbolically for $f(t)$.
4. Apply the Laplace transform to obtain the final expression for $u(x)$.

This process naturally yields symbolic series expansions involving terms such as $x^{\alpha}\ln(x)^{n}$, which are challenging to obtain using standard numerical or semi-analytical methods.

### 4.2. Algorithmic Steps

The procedure can be formally summarized as follows.

**Algorithm 1: RILT Method**

Input: IVP in the form $x^{\alpha(x)}D^{\alpha(x)}\big(u(x)\big) = A(u(x), x)$.

Output: Symbolic approximation of the solution $u(x)$.

**Step 1: Assume solution form:**

Let $u(x) = \mathcal{L}\{f(t)\} = \int_0^{\infty} f(t)x^{t}dt$.

**Step 2: Apply inverse Laplace transform:**

Invert both sides of the equation to express the problem in terms of $f(t)$:

$$\mathcal{L}^{-1}\{x^{\alpha(x)}D^{\alpha(x)}u(x)\} = \mathcal{L}^{-1}\{A(u(x), x)\}$$

$$\Rightarrow f(t) = \Theta\big(t - \alpha(x)\big)\frac{\Gamma(t - \alpha(x) + 1)}{\Gamma(t + 1)}\mathcal{L}^{-1}\{A(u(x), x)\}.$$

**Step 3: Imposing initial conditions:**

$$f(t) = \sum_{k=0}^{n-1}\frac{c_k}{k!}\delta(t - k) + \Theta\big(t - \alpha(x)\big)\frac{\Gamma(t - \alpha(x) + 1)}{\Gamma(t + 1)}\mathcal{L}^{-1}\{A(u(x), x)\}.$$

**Step 4: Applying Laplace transform:**

$$u(x) = \mathcal{L}\{f(t)\} = \sum_{k=0}^{n-1}\frac{c_k}{k!}x^{k} + \mathcal{L}\left\{\Theta\big(t - \alpha(x)\big)\frac{\Gamma(t - \alpha(x) + 1)}{\Gamma(t + 1)}\mathcal{L}^{-1}\{A(u(x), x)\}\right\}.$$

The second term in the right hand side is found by using Lemmas 2 and 3 after expanding $A(u(x), x)$ as logarithmic-power series if required.

**Step 5: Recursive solving for an arbitrary accuracy:**

$$u_{j+1}(x) = \sum_{k=0}^{n-1}\frac{c_k}{k!}x^{k} + \mathcal{L}\left\{\Theta\big(t - \alpha(x)\big)\frac{\Gamma(t - \alpha(x) + 1)}{\Gamma(t + 1)}\mathcal{L}^{-1}\{A\big(u_j(x), x\big)\}\right\}$$

$$\Rightarrow u_N(x) = \sum_{k=0}^{N}\sum_{r=0}^{N} c_{k,r} x^{\frac{k}{q}} \ln(x)^r, \qquad q \in \mathbb{Z}^+.$$

### 4.3. Convergence of RILT

To establish the convergence of the symbolic series produced by RILT method, we interpret the solution $u(x)$ as a Laplace-type integral of the form:

$$u(x) = \int_0^\infty f(t) x^t dt, \qquad 0 < x < 1,$$

where $f(t)$ is the auxiliary function obtained via the recursive process described in Section 4.2. The approximation generated by the method can then be viewed as the truncated integral:

$$u_N(x) = \int_0^N f(t) x^t dt,$$

where $N \in \mathbb{R}^+$ denotes the truncation level.

The truncation error is defined as:

$$E_N(x) = u(x) - u_N(x) = \int_N^\infty f(t) x^t dt.$$

From the necessary condition of existence of Laplace transform, there exist constants $M, \beta > 0$ such that $|f(t)| \leq Me^{-\beta}$ for all $t \geq 0$.

Under this condition, and for $x \in (0, e^\beta)$, the integrand satisfies:

$$|f(t)x^t| \leq M\left(xe^{-\beta}\right)^t.$$

Therefore, the tail integral converges and admits the following bound:

$$|E_N(x)| = \left|\int_N^\infty f(t)x^t dt\right| \leq \int_N^\infty |f(t)x^t| dt = \int_N^\infty M\left(xe^{-\beta}\right)^t dt = -\frac{M\left(xe^{-\beta}\right)^N}{\ln(e^{-\beta}x)}.$$

For $x < e^\beta$:

$$\lim_{N\to\infty} |E_N(x)| \leq 0$$
$$\Rightarrow u_N(x) \to u(x) \text{ as } N \to \infty.$$

### 4.4. Strengths and Applicability

The RILT method is particularly suited for:

- Equations with fractional or variable-order derivatives.
- Problems exhibiting singularities or requiring non-integer power-log terms.
- Equations where symbolic dependence on parameters is essential.
- High-precision solutions with arbitrary accuracy (due to symbolic recursive structure).

- Systems with delay or memory terms, when kernel forms are Laplace-invertible.

The method overcomes limitations typically encountered in generating Adomian polynomials or solving nonlinear algebraic systems at each step. Its symbolic nature allows for complete control over accuracy and form, and it can be implemented using standard computer algebra systems (e.g., Maple, Mathematica, Sympy).

## 5. Examples

In this section, we demonstrate the accuracy, generality, and symbolic capability of the proposed RILT method through a series of benchmark problems. The selected examples include linear and nonlinear differential equations, equations with singularities, fractional and variable-order models, and problems involving delay or integral terms.

Each example is solved symbolically using the RILT method, and where appropriate, the results are compared with known exact solutions or standard numerical methods such as the classical fourth-order Runge–Kutta (RK4) scheme. In all cases, the symbolic nature of the method allows for direct manipulation of parameters and arbitrary expansion order without relying on numerical discretization or approximation.

**Example 1: Linear First-Order Differential Equation**

The simplest case that doesn't require recursive solving is the linear DE

$$x^2 y''(x) - 3xy'(x) - 3y(x) = \frac{1}{\ln(x)^2}.$$

Using Lemma 1, the equation is written as

$$\mathcal{L}\{t(t-1)f(t) - 3tf(t) - 3f(t)\} = \mathcal{L}\{t\},$$

which simplifies to

$$(t^2 - 4t - 3)f(t) = t.$$

The exact particular solution is then obtained by

$$y(x) = \mathcal{L}\{f(t)\} = \mathcal{L}\left\{\frac{t}{t^2 - 4t - 3}\right\} = \int_0^\infty \frac{t}{t^2 - 4t - 3} x^t dt,$$

$$y(x) = \frac{x^{2-\sqrt{7}}\left(\left(\sqrt{7}+2\right)x^{2\sqrt{7}}\,\mathrm{Ei}\left(-\left(\left(\sqrt{7}+2\right)\ln(x)\right)\right) + \left(\sqrt{7}-2\right)\mathrm{Ei}\left(\left(\sqrt{7}-2\right)\ln(x)\right)\right)}{2\sqrt{7}},$$

where $Ei(x)$ is the exponential integral function defined as

$$Ei(x) = \int_{-\infty}^{x} \frac{e^t}{t} dt.$$

**Example 2: Nonlinear Riccati Differential Equation**

$$y'(x) - y^2(x) = 1.$$

Step 1: Write the differential equation in the standard form (3):

$$xy'(x) = x + xy^2(x).$$

Step 2: Apply the inverse Laplace transform of the equation:

$$tf(t) = \mathcal{L}^{-1}\{x + xy^2(x)\} \Rightarrow f(t) = \frac{\mathcal{L}^{-1}\{x + xy^2(x)\}}{t}.$$

Step 3: Impose the initial condition $y(0) = c$:

$$f(t) = c\delta(t) + \frac{\mathcal{L}^{-1}\{x + xy^2(x)\}}{t}$$

Step 4: Find the final recurrence relation:

$$y_{j+1}(x) = c + \mathcal{L}\left\{\frac{\mathcal{L}^{-1}\{x + xy_j^2(x)\}}{t}\right\}.$$

Step 5: Solve the recurrence relation up to the desired degree:

$$y_0(x) = c,$$

$$y_1(x) = c + \mathcal{L}\left\{\frac{\mathcal{L}^{-1}\{x + y_0^2(x)x\}}{t}\right\} = c + \mathcal{L}\left\{\frac{(1 + c^2)\delta(t-1)}{t}\right\} = c + (1 + c^2)x,$$

$$y_2(x) = c + \mathcal{L}\left\{\frac{\mathcal{L}^{-1}\{x + y_1^2(x)x\}}{t}\right\} = c + (1 + c^2)x + c(1 + c^2)x^2,$$

$$y_3(x) = c + \mathcal{L}\left\{\frac{\mathcal{L}^{-1}\{x + y_2^2(x)x\}}{t}\right\} = c + (1 + c^2)x + c(1 + c^2)x^2 + \left(\frac{1}{3} + \frac{4}{3}c^2 + c^4\right)x^3,$$

where $y_n(x)$ is the solution up to $x^n$.

The general solution acquired is the Maclaurin series of $y(x) = \tan(x + \tan^{-1}(c))$.

**Example 3: Nonlinear Riccati Differential Equation**

$$y'(x) - y^2(x) = x^{\sin(x)-1}\left(-x^{\sin(x)+1} + \sin(x) + x\ln(x)\cos(x)\right).$$

Step 1: Write the differential equation in the standard form (3):

$$xy'(x) = x^{\sin(x)}\left(-x^{\sin(x)+1} + \sin(x) + x\ln(x)\cos(x)\right) + xy^2(x).$$

Step 2: Apply the inverse Laplace transform of the equation:

$$tf(t) = \mathcal{L}^{-1}\{g(x) + xy^2(x)\} \Rightarrow f(t) = \frac{\mathcal{L}^{-1}\{g(x) + xy^2(x)\}}{t},$$

where $g(x) = x^{\sin(x)}\left(-x^{\sin(x)+1} + \sin(x) + x\ln(x)\cos(x)\right)$.

Step 3: Impose the initial condition $y(0) = c$:

$$f(t) = c\delta(t) + \frac{\mathcal{L}^{-1}\{g(x) + xy^2(x)\}}{t}$$

Step 4: Find the final recurrence relation:

$$y_{j+1}(x) = c + \mathcal{L}\left\{\frac{\mathcal{L}^{-1}\{g(x) + xy_j^2(x)\}}{t}\right\}.$$

Step 5: Solve the recurrence relation, after expanding $g(x)$ as logarithmic-power series up to the desired degree:

$$y_0(x) = c,$$

$$y_1(x) = c + \mathcal{L}\left\{\frac{\mathcal{L}^{-1}\{g(x) + xy_0^2(x)\}}{t}\right\} = c + \mathcal{L}\left\{\frac{\mathcal{L}^{-1}\{x(c^2 + ln(x))\}}{t}\right\}$$
$$= c + \mathcal{L}\left\{\frac{c^2\delta(t-1) + \delta'(t-1)}{t}\right\} = c + (c^2 - 1)x + xln(x),$$

$$y_2(x) = c + \mathcal{L}\left\{\frac{\mathcal{L}^{-1}\{g(x) + xy_1^2(x)\}}{t}\right\}$$
$$= c^3x^2 + c^2x - \frac{3cx^2}{2} + x\left((c-1)x + 1\right)\ln(x) + c + \frac{1}{2}x^2\ln(x)^2 + \frac{1}{2}(x-2)x,$$

$$y_3(x) = c + \mathcal{L}\left\{\frac{\mathcal{L}^{-1}\{g(x) + xy_2^2(x)\}}{t}\right\}$$
$$= x^3\left(c^4 + \frac{1}{18}(24c^2 - 16c - 11)\ln(x) - \frac{19c^2}{9} + \frac{1}{3}(c-1)\ln(x)^2 + \frac{17c}{27} + \frac{\ln(x)^3}{6} + \frac{13}{27}\right) + x^2\left(c^3 + (c-1)\ln(x) - \frac{3c}{2} + \frac{\ln(x)^2}{2} + \frac{1}{2}\right) + x(c^2 + \ln(x) - 1) + c.$$

If the initial condition $y(0) = 1$ is applied, the logarithmic-power series of the exact solution $y(x) = x^{\sin(x)}$ could be found up to any desired accuracy.

**Example 4: Coupled Asymmetric Van der Pol Oscillator** [7]

Consider the two-dimensional coupled asymmetric Van der Pol oscillator:

$$x''(t) + \omega^2 x(t) + \epsilon(\alpha x(t)^2 + \beta y(t)^2 - 1)x'(t) = 0,$$
$$y''(t) + \omega^2 y(t) - \epsilon(\alpha x(t)^2 + \beta y(t)^2 - 1)y'(t) = 0,$$

with initial conditions:

$$x(0) = x_0, \qquad x'(0) = x_1, \qquad y(0) = y_0, \qquad y'(0) = y_1.$$

In 0.24 seconds, the proposed method provides a third-degree symbolic approximation to the solution.

Approximate symbolic solution (up to degree 3):

$$\begin{aligned} x(t) = x_0 + tx_1 - \frac{1}{2}t^2\big(\omega^2 x_0 + \alpha\epsilon x_0^2 x_1 + \epsilon x_1(-1 + \beta y_0^2)\big) \\ + \frac{1}{6}t^3\big(\alpha\epsilon\omega^2 x_0^3 + \alpha^2\epsilon^2 x_0^4 x_1 - 2\alpha\epsilon x_0 x_1^2 + \epsilon\omega^2 x_0(-1 + \beta y_0^2) \\ + 2\alpha\epsilon^2 x_0^2 x_1(-1 + \beta y_0^2) + x_1(\epsilon^2 - \omega^2 - 2\beta\epsilon^2 y_0^2 + \beta^2\epsilon^2 y_0^4 - 2\beta\epsilon y_0 y_1)\big), \end{aligned}$$

$$\begin{aligned} y(t) = y_0 + ty_1 + \frac{1}{2}t^2(-\omega^2 y_0 + \epsilon(-1 + \alpha x_0^2)y_1 + \beta\epsilon y_0^2 y_1) \\ + \frac{1}{6}t^3(\epsilon\omega^2(1 - \alpha x_0^2)y_0 - \beta\epsilon\omega^2 y_0^3 \\ + (\epsilon^2 - \omega^2 - 2\alpha\epsilon^2 x_0^2 + \alpha^2\epsilon^2 x_0^4 + 2\alpha\epsilon x_0 x_1)y_1 + 2\beta\epsilon^2(-1 + \alpha x_0^2)y_0^2 y_1 \\ + \beta^2\epsilon^2 y_0^4 y_1 + 2\beta\epsilon y_0 y_1^2). \end{aligned}$$

For the values $\alpha = 1.546, \beta = 1.551, \omega = \epsilon = 1, x(0) = y(0) = 1, x'(0) = y'(0) = 0$, the approximate solution up to degree 20 is found in 1.2 seconds as following

$$\begin{aligned} x(t) \\ = 1 - \frac{t^2}{2} + \frac{699t^3}{2000} - \frac{3397409t^4}{24000000} - \frac{4518211109t^5}{40000000000} + \frac{109670503086719t^6}{720000000000000} - \cdots \\ + \frac{97293050284120683172229514052034920532533219807693041285194013323053t^{20}}{130101711667200000000000000000000000000000000000000000000000000000000000}, \end{aligned}$$

$$\begin{aligned} y(t) \\ = 1 - \frac{t^2}{2} - \frac{699t^3}{2000} - \frac{3397409t^4}{24000000} + \frac{4518211109t^5}{40000000000} + \frac{109838263086719t^6}{720000000000000} + \cdots \\ + \frac{18377460773033817205214713729400971331841903136360609292620893531410911t^{20}}{243290200817664000000000000000000000000000000000000000000000000000000000}. \end{aligned}$$

Due to the finite radius of convergence of the Taylor series and the need to obtain a solution over a larger time domain, a successive Taylor expansion is performed using a small time step $\Delta t = \frac{1}{10}$. The corresponding solution profiles are shown in Figures 1 and 2, while the phase portrait is presented in Figure 3.

The results show an excellent agreement with those reported by Erturk et al. on the same problem.

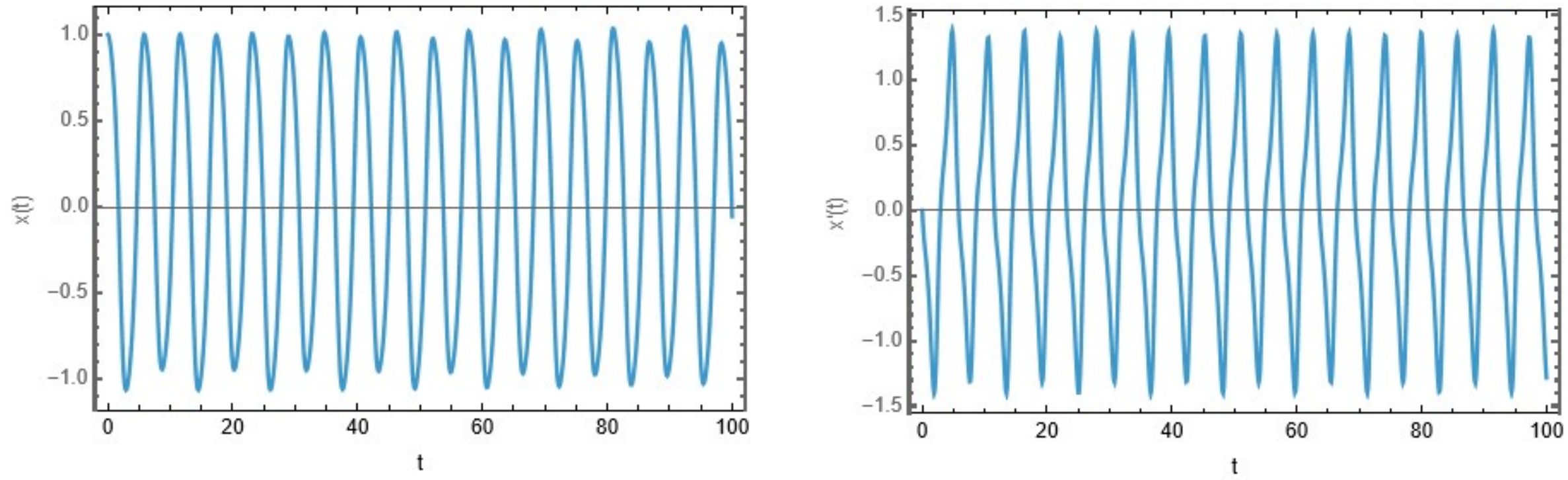

**Fig 1** The solution $x(t)$ of Example 4 and $x'(t)$ for $0 \leq t \leq 100$

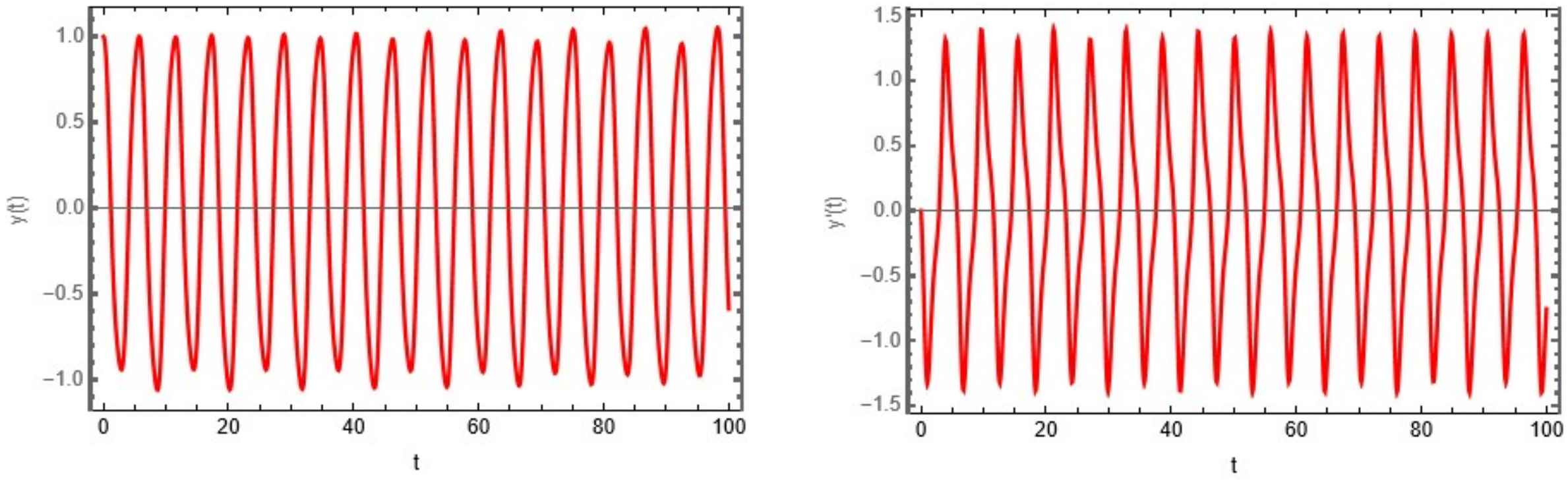

**Fig 2** The solution $y(t)$ of Example 4 and $y'(t)$ for $0 \leq t \leq 100$

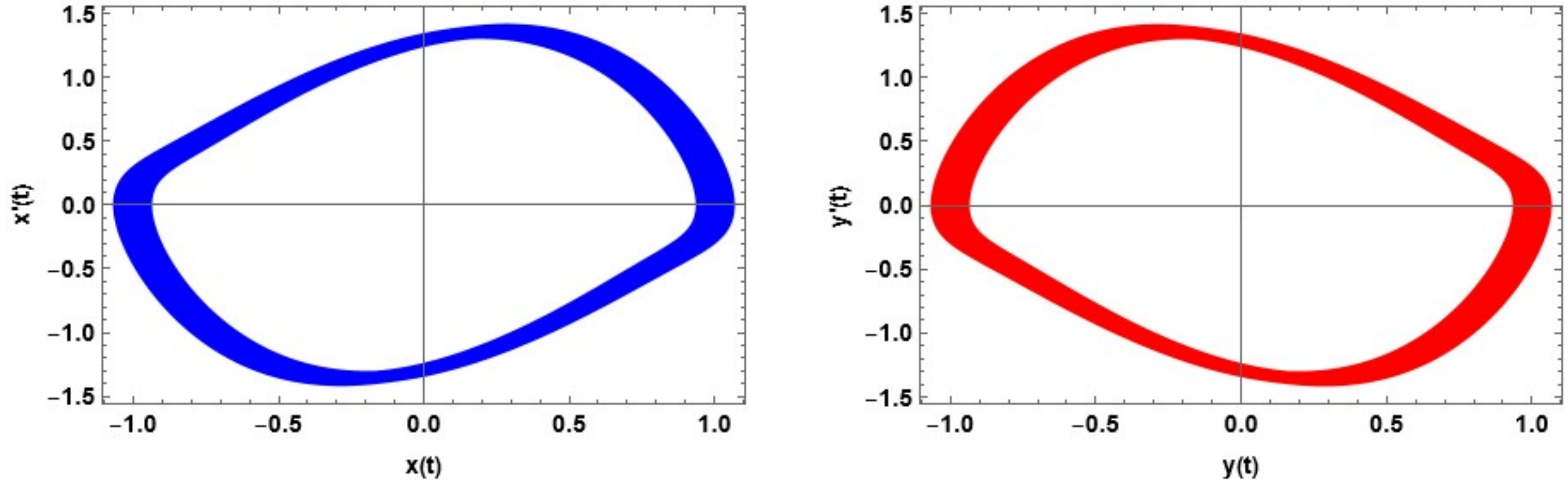

**Fig 3** The phase diagrams of Example 4 for $0 \leq t \leq 300$

**Example 5: Inverted Rössler System** [8]

Consider the inverted Rössler system

$$\dot{X}(t) = -Y(t) - X(t)Y(t)Z(t),$$

$$\dot{Y}(t) = X(t) + AY(t),$$

$$\dot{Z}(t) = BX(t) + Z(t)(X(t)^2 - C),$$

with initial conditions:

$$X(0) = X_0, \qquad Y(0) = Y_0, \qquad Z(0) = Z_0,$$

where $A, B$ and $C$ are the control parameters.

In 0.44 seconds, the proposed method computes a third-degree symbolic approximation of the solution.

Approximate symbolic solution (up to degree 3):

$$\begin{aligned} X(t) = X_0 &- t(Y_0 + X_0 Y_0 Z_0) \\ &- \frac{1}{2}t^2 (X_0 + AY_0 + X_0^2(BY_0 + Z_0) + X_0Y_0Z_0(A - C) + X_0^3 Y_0 Z_0 - Y_0^2 Z_0 \\ &- X_0 Y_0^2 Z_0^2) \\ &+ \frac{1}{6}t^3 \Big(-AX_0 - 2BX_0^3 + Y_0 - A^2Y_0 - 2ABX_0^2Y_0 + BCX_0^2Y_0 + 3BX_0Y_0^2 - AX_0^2Z_0 \\ &+ 2CX_0^2Z_0 + 5X_0Y_0Z_0 - A^2X_0Y_0Z_0 + 2ACX_0Y_0Z_0 - C^2X_0Y_0Z_0 - 2AX_0^3Y_0Z_0 \\ &+ 2CX_0^3Y_0Z_0 - X_0^5Y_0Z_0 + 3AY_0^2Z_0 - 2CY_0^2Z_0 + 4(1+B)X_0^2Y_0^2Z_0 + 4X_0^2Y_0Z_0^2 \\ &+ 3AX_0Y_0^2Z_0^2 - 3CX_0Y_0^2Z_0^2 + 5X_0^3Y_0^2Z_0^2 - Y_0^3Z_0^2 - X_0Y_0^3Z_0^3 - X_0^4(BY_0 + 2Z_0)\Big), \end{aligned}$$

$$\begin{aligned} Y(t) = Y_0 &+ t(X_0 + AY_0) + \frac{1}{2}t^2(AX_0 - Y_0 + A^2Y_0 - X_0Y_0Z_0) \\ &+ \frac{1}{6}t^3\big(-X_0 + A^2X_0 - 2AY_0 + A^3Y_0 - 2AX_0Y_0Z_0 + CX_0Y_0Z_0 - X_0^3Y_0Z_0 + Y_0^2Z_0 \\ &+ X_0Y_0^2Z_0^2 - X_0^2(BY_0 + Z_0)\big), \end{aligned}$$

$$\begin{aligned} Z(t) = Z_0 &+ t(BX_0 + (X_0^2 - C)Z_0) \\ &+ \frac{t^2}{2}(B(X_0^3 - Y_0 - CX_0) + (C^2 - 2CX_0^2 + X_0^4)Z_0 - (2+B)X_0Y_0Z_0 - 2X_0^2Y_0Z_0^2) \\ &+ \frac{t^3}{6}((C^2-1)BX_0 - 2BCX_0^3 + BX_0^5 - ABY_0 + BCY_0 - B(5+B)X_0^2Y_0 - C^3Z_0 \\ &+ (-2-B)X_0^2Z_0 + 3C^2X_0^2Z_0 - 3CX_0^4Z_0 + X_0^6Z_0 - 2AX_0Y_0Z_0 - ABX_0Y_0Z_0 \\ &+ 6CX_0Y_0Z_0 + 2BCX_0Y_0Z_0 - 2(3+4B)X_0^3Y_0Z_0 + (2+B)Y_0^2Z_0 - 2X_0^3Z_0^2 \\ &- 2AX_0^2Y_0Z_0^2 + 8CX_0^2Y_0Z_0^2 - 8X_0^4Y_0Z_0^2 + (6+B)X_0Y_0^2Z_0^2 + 4X_0^2Y_0^2Z_0^3). \end{aligned}$$

For the particular values $A = B = \frac{1}{5}$ and $C = \frac{11}{2}$, and for the initial conditions (1,0,0), the approximate solution up tp degree 20 is found in 1.17 seconds as following

$$\begin{aligned} X(t) = 1 - \frac{t^2}{2} - \frac{t^3}{10} + \frac{59t^4}{400} - \frac{259t^5}{3000} + \frac{254477t^6}{3600000} - \frac{908857t^7}{15750000} + \cdots \\ + \frac{112443221601004180761044729272 51t^{20}}{477039609446400000000000000000}, \end{aligned}$$

$$Y(t) = t + \frac{t^2}{10} - \frac{4t^3}{25} - \frac{33t^4}{1000} + \frac{1409t^5}{50000} - \frac{60523t^6}{4500000} + \frac{6119833t^7}{630000000} - \cdots - \frac{76742393900037702819664130890069t^{20}}{50685458503680000000000000000000000},$$

$$Z(t) = \frac{t}{5} - \frac{9t^2}{20} + \frac{77t^3}{120} - \frac{1243t^4}{1600} + \frac{62967t^5}{80000} - \frac{1920877t^6}{2880000} + \frac{490214341t^7}{1008000000} - \cdots + \frac{20258110899918487169925963083629t^{20}}{22248760934400000000000000000000},$$

After successive Taylor expansion of the solution with $\Delta t = \frac{1}{10}$, the behavior of the system for the domain $0 \le t \le 800$ could be seen in Figure 4. Table 1 shows the solution obtained for $0 \le t \le 2$ and Table 2 shows the comparison of the acquired solution and the solution obtained by $RK4$ and $NCFD$ which proves the reliability of the method.

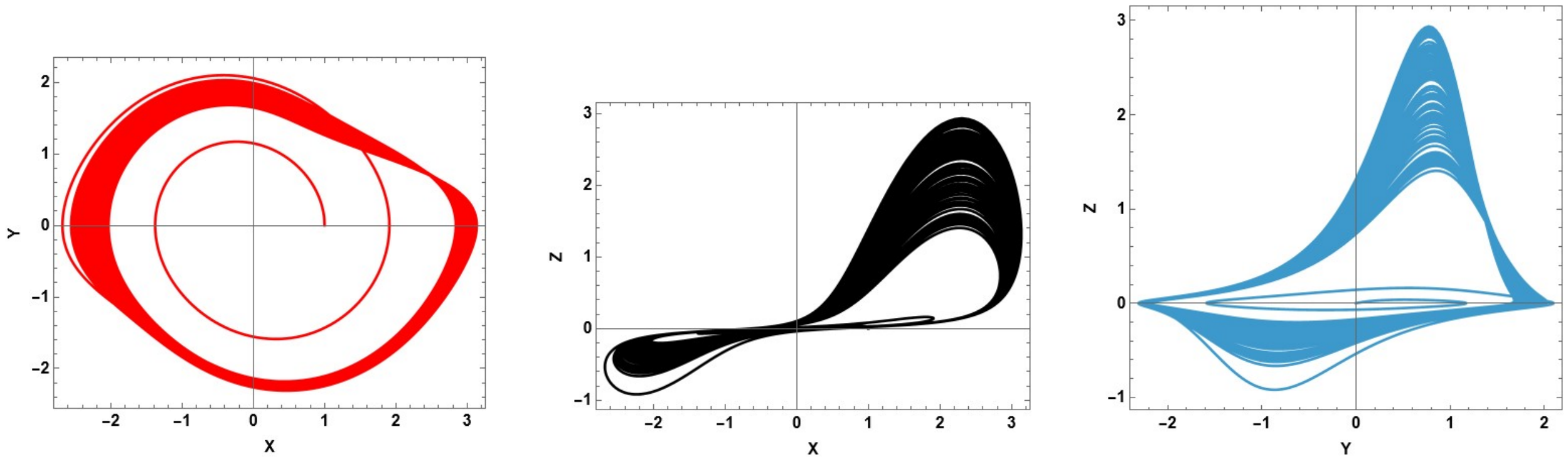

**Fig 4** Views of Example 5 for $0 \le t \le 800$

**Example 6: Nonlinear Singular IVP** [9]

Consider the nonlinear singular IVP

$$\frac{1}{x^2}\frac{d}{dx}\left(x^2\frac{du}{dx}\right) = -u^5, \qquad u(0) = 1, \qquad u'(0) = 0, \qquad 0 \le x \le 1,$$

With the exact solution $u(x) = \frac{1}{\sqrt{1+\frac{x^2}{3}}}$.

In 0.1 seconds, the proposed method found the tenth-degree symbolic approximate solution:

$$u(x) = c_0 + c_1 x - \frac{{c_0}^5 x^2}{6} - \frac{5}{12}{c_0}^4 c_1 x^3 + \left(\frac{{c_0}^9}{24} - \frac{{c_0}^3 {c_1}^2}{2}\right)x^4 + \left(\frac{13{c_0}^8 c_1}{72} - \frac{{c_0}^2 {c_1}^3}{3}\right)x^5 + \cdots + \left(-\frac{7c_0^{21}}{6912} + \frac{27289c_0^{15}c_1^2}{266112} - \frac{23c_0^9c_1^4}{42} + \frac{5c_0^3c_1^6}{42}\right)x^{10},$$

which after imposing the initial conditions, the particular solution is found. Table 3 shows the comparison of the maximum absolute errors of the solutions obtained by $RILT$, $RK4$ and one step method with the time of $RILT$ at different number of terms.

**Table 1** Solution of Example 5 for $0 \leq t \leq 2$

| $t$ | $X(t)$ | $Y(t)$ | $Z(t)$ |
|---|---|---|---|
| 0.0 | 1.0 | 0.0 | 0.0 |
| 0.1 | 0.99491395194157500 | 0.10083696925599112 | 0.016071228793570442 |
| 0.2 | 0.97941224516387910 | 0.20267546519417630 | 0.026105031359680817 |
| 0.3 | 0.95332597299918400 | 0.30447310832863540 | 0.032047763218575630 |
| 0.4 | 0.91670693099270040 | 0.40516089369098890 | 0.035138346329579490 |
| 0.5 | 0.86978832412060770 | 0.50366287038371230 | 0.036188494412550530 |
| 0.6 | 0.81295398823078350 | 0.59891272493269290 | 0.035750651912722635 |
| 0.7 | 0.74671229563874520 | 0.68986784630084470 | 0.034217451430580580 |
| 0.8 | 0.67167338885942730 | 0.77552121131898980 | 0.031879438313269184 |
| 0.9 | 0.58852965178016630 | 0.85491137168328700 | 0.028957765895684242 |
| 1.0 | 0.49803970949498977 | 0.92713084542528380 | 0.025622449573619780 |
| 1.1 | 0.40101608921914267 | 0.99133324640766120 | 0.022002838716397860 |
| 1.2 | 0.29831630087102230 | 1.04673948730368840 | 0.018194314566499100 |
| 1.3 | 0.19083675367179054 | 1.09264335568237580 | 0.014263414479875391 |
| 1.4 | 0.07950874469200490 | 1.12841669924595430 | 0.010252404665730061 |
| 1.5 | -0.034704242730112675 | 1.15351438256192230 | 0.006183635464934206 |
| 1.6 | -0.150808558178109340 | 1.16747910990088230 | 0.002063686896263901 |
| 1.7 | -0.267782155036066400 | 1.16994615686458640 | -0.0021127774426596743 |
| 1.8 | -0.384576359200522970 | 1.16064802022197780 | -0.0063594913252755240 |
| 1.9 | -0.500117740208297100 | 1.13941897832357180 | -0.0106953567552678180 |
| 2.0 | -0.613310113496490800 | 1.10619954815831900 | -0.0151413694183409330 |

**Table 2** Comparison of the solution of RILT with $RK4$ and $NCFD$ for $t = 2, (X, Y, Z) = (1,0,0)$

| | $X(t)$ | $Y(t)$ | $Z(t)$ |
|---|---|---|---|
| $RK4\left(\Delta t = \frac{1}{10000}\right)$ | -0.611628121823525 | 1.112389089856173 | -0.008030733489998 |
| $NCFD\left(\Delta t = \frac{1}{10240}\right)$ | -0.603158345062954 | 1.053411896544050 | -0.007933575350318 |
| $RILT\left(\Delta t = \frac{1}{10}\right)$ | -0.6133101134964908 | 1.106199548158319 | -0.015141369418340 |
| $RILT\left(\Delta t = \frac{1}{10000}\right)$ | -0.6133101134962827 | 1.106199548158433 | -0.015141369418331 |

**Table 3** Comparison of the maximum absolute error of $RILT$, $RK4$ and $one\ step\ method$

| $N$ | $RK4$ | $One\ step\ method$ | $RILT$ | $Time\ (s)$ |
|---|---|---|---|---|
| 16 | $6.76 \times 10^{-6}$ | $1.56 \times 10^{-5}$ | $7.19 \times 10^{-7}$ | 0.008 |
| 64 | $1.05 \times 10^{-7}$ | $9.29 \times 10^{-7}$ | $1.32 \times 10^{-17}$ | 0.231 |
| 256 | $1.65 \times 10^{-9}$ | $6.11 \times 10^{-8}$ | $1.05 \times 10^{-63}$ | 9.9 |
| 2048 | $3.23 \times 10^{-12}$ | | | |
| 65536 | | $9.83 \times 10^{-13}$ | | |

**Example 7: Nonlinear Fractional Order Riccati Differential Equation** [10]

Consider the nonlinear fractional order Riccati differential equation

$$D^{\alpha}x(t) + x^2(t) - 1 = 0, \qquad x(0) = 0, \qquad 0 \le t \le 1,$$

with the exact solution $x(t) = \tanh(t)$ at $\alpha = 1$.

In 2.1 seconds, the proposed method found the symbolic approximate solution:

$$\begin{aligned} x(t) = x_0 &+ \frac{1 - x_0^2}{\Gamma(1+\alpha)} t^{\alpha} + \frac{2x_0(-1 + x_0^2)}{\Gamma(1+2\alpha)} t^{2\alpha} \\ &- \frac{(-1 + x_0^2)\left(4x_0^2\Gamma(1+\alpha)^2 + (-1 + x_0^2)\Gamma(1+2\alpha)\right)}{\Gamma(1+\alpha)^2\Gamma(1+3\alpha)} t^{3\alpha} \\ &+ \frac{2x_0(-1 + x_0^2)}{\Gamma(1+\alpha)^2\Gamma(1+2\alpha)\Gamma(1+4\alpha)} \left(4x_0^2\Gamma(1+\alpha)^2\Gamma(1+2\alpha)\right. \\ &+ \left.(-1 + x_0^2)\Gamma(1+2\alpha)^2 + 2(-1 + x_0^2)\Gamma(1+\alpha)\Gamma(1+3\alpha)\right) t^{4\alpha}, \end{aligned}$$

which is a function of the fractional order $\alpha$.

For $\alpha = 1$ and $x_0 = 0$, the approximate solution and corresponding absolute error and time are shown in Table 4.

**Table 4** Maximum absolute error of $RILT$ and time at different number of terms for $\alpha = 1, x_0 = 0$

| $N$ | $\lvert Error \rvert_{Max}$ | $Time\ (s)$ |
|---|---|---|
| 20 | $6.89 \times 10^{-5}$ | 0.02 |
| 60 | $9.85 \times 10^{-13}$ | 0.17 |
| 200 | $1.11 \times 10^{-16}$ | 4.61 |

The maximum absolute errors obtained by the Mittag-Leffler-Galerkin method and the algorithm proposed by Ezz-Eldien et al.[11] are in order of $10^{-1}$ .

**Example 8: Linear Fractional Volterra Integro-Differential Equation** [12]

Consider the linear fractional integro-differential equation with the weakly singular kernel

$$D^{\frac{1}{3}}u(x) = \int_0^x \frac{D^{\frac{1}{3}}u(z)}{(x-z)^{\frac{1}{2}}} dz + \frac{6}{\Gamma\left(\frac{11}{3}\right)} x^{\frac{8}{3}} - \frac{6\sqrt{\pi}}{\Gamma\left(\frac{25}{6}\right)} x^{\frac{19}{6}}, \qquad 0 \le x \le 1,$$

with the initial condition $u(0) = 0$ and the exact solution $u(x) = x^3$.

The exact solution is found by RILT in $0.1$ seconds, while the maximum absolute error obtained by SSKCPs collocation method is $3.8206 \times 10^{-6}$ in $1.576$ seconds.

**Example 9: Nonlinear Fractional Volterra-Fredholm Integro-Differential Equation** [12]

Consider the nonlinear fractional integro-differential equation with the weakly singular kernel

$$D^{\alpha}u(x) = \int_0^x \frac{u^2(z)}{(x-z)^{\frac{1}{2}}}dz + \int_0^1 xzu^2(z)dz + 3x^2 - \frac{\Gamma(7)\Gamma\left(\frac{1}{2}\right)}{\Gamma\left(\frac{15}{2}\right)}x^{\frac{13}{2}} - \frac{x}{8}, \qquad 0 \le x \le 1,$$

with the initial condition $u(0) = 0$ and the exact solution is $u(x) = x^3$ at $\alpha = 1$.

The Fredholm term in the equation is introduced in the algorithm as $cx$ where $c$ is a constant to be found later. The problem now becomes

$$D^{\alpha}u(x) = \int_0^x \frac{u^2(z)}{(x-z)^{\frac{1}{2}}}dz + cx + 3x^2 - \frac{\Gamma(7)\Gamma\left(\frac{1}{2}\right)}{\Gamma\left(\frac{15}{2}\right)}x^{\frac{13}{2}} - \frac{x}{8}, \qquad 0 \le x \le 1,$$

which is solved by the proposed method to get the solution

$$u(x) = x^3 + \left(\frac{c}{2} - \frac{1}{16}\right)x^2.$$

The next step is to find the value of $c$ which could be done by plugging the solution in the original equation to get $c = \frac{1}{8}$ and hence the exact solution is found in $0.72$ seconds, while the maximum absolute error obtained by SSKCPs collocation method is $4.8129 \times 10^{-14}$ in $95.784$ seconds.

**Example 10: Nonlinear Variable-Order Fractional Differential Equation** [13]

Consider the nonlinear variable order fractional differential equation

$$D^{\alpha(x)}u(x) + \sin(x)\, u^2(x) = 6\frac{x^{3-\alpha(x)}}{\Gamma\big(4-\alpha(x)\big)} + x^6 \sin(x), \qquad 0 \le x \le 1,$$

with the initial condition $u(0) = 0$, $\alpha(x) = 1 - \frac{e^{-x}}{2}$ and the exact solution $u(x) = x^3$.

The exact solution is found by RILT in $0.008$ seconds, while the maximum absolute error of the approximate solution obtained by SMLRM-CPs is $2.2696 \times 10^{-12}$.

**Example 11: Delay Variable-Order Fractional Differential Equation** [13]

Consider the delay variable order fractional differential equation

$$D^{\alpha(x)}u(x) + e^x u'(x) + 2u(x) + 8u(x-1) = g(x), \qquad 0 \le x \le 1,$$

with the boundary conditions $u(0) = 0, u(1) = 9$, $\alpha(x) = \frac{6+\cos(x)}{4}$.

The non-homogeneous part is $g(x) = (4+2x)e^x + 16 + 24x + 10x^2 + 2\frac{x^{2-\alpha(x)}}{\Gamma(3-\alpha(x))}$, and the exact solution is $u(x) = x^2 + 4x + 4$.

In 0.015 seconds, RILT found the solution $u(x) = x^2 + cx + 4$ which after imposing the boundary condition, the exact solution is found. The maximum absolute error of the approximate solution obtained by SMLRM-CPs is $1.3257 \times 10^{-13}$.

**Example 12: Delay Fractional Order Differential Equation** [14]

Consider the delay fractional order differential equation

$$D^{\rho}w(z) = w(z) + \frac{1}{10}w\left(\frac{1}{10}z\right) + g(z), \qquad 0 \le z \le 1, \qquad 0 < \rho \le 1,$$

with the initial condition $w(0) = 1$, $g(z) = \frac{2\rho}{\Gamma(3-\rho)}z^{2-\rho} - \frac{11}{10} - \rho z^2 - \frac{\rho}{1000}z^2$ and the exact solution $w(z) = \rho z^2 + 1$.

In 0.03 seconds, RILT found the symbolic exact solution $w(z) = \rho z^2 + c$. The maximum absolute error of the approximate solution obtained by GNPs matrix collocation method is in order of $10^{-16}$ for $\rho \neq 1$.

**Example 13: Nonlinear Delay Variable-Order Fractional Differential Equation**

Consider the following nonlinear delay variable order fractional differential equation

$$D^{\alpha(x)}y(x) + \left(y\left(\frac{x-1}{3}\right)\right)^2 + y(y(x)) = g(x),$$

where

$$g(x) = \frac{1}{6561}\left(\begin{array}{c} \frac{157464x^{\frac{e^{-x}}{2}+3}}{\Gamma\left(4+\frac{e^{-x}}{2}\right)} - \frac{32805x^{\frac{e^{-x}}{2}}}{\Gamma\left(1+\frac{e^{-x}}{2}\right)} + 6561x^{16} - 131220x^{13} + 52488x^{12} \\ +984150x^{10} - 787320x^9 + 157465x^8 - 3280508x^7 + 3936628x^6 \\ -1574966x^5 + 4279516x^4 - 6565052x^3 + 3959497x^2 - 968579x + 141292 \end{array}\right),$$

with the initial condition $y(0) = 2$ and the exact solution $y(x) = x^4 - 5x + 2$ at the variable order $\alpha(x) = 1 - \frac{e^{-x}}{2}$.

The exact solution is found by RILT in 1.46 seconds.

## 6. Conclusion

In this work, we introduced a novel semi-analytical technique called the Recursive Inverse Laplace Transform (RILT) for solving a wide class of initial value problems, including differential, integro-differential, fractional, and delay equations. Unlike conventional Laplace-based methods that begin with forward transformation, the proposed method reverses the procedure by first applying the inverse Laplace transform to the equation. This enables symbolic construction of the solution in the Laplace time domain, followed by a final Laplace transform to obtain the desired approximate solution.

The method is capable of generating symbolic solutions with arbitrary accuracy, often in the form of generalized logarithmic-power series. It naturally accommodates challenging features such as variable-order derivatives, singularities, memory effects, and time delays, which are typically difficult to address using standard numerical or semi-analytical techniques.

A wide variety of examples were presented to validate the method, ranging from simple linear equations to chaotic systems, nonlinear oscillators, fractional and variable-order models, and coupled systems. In each case, the symbolic series solution obtained by the RILT method showed excellent agreement with known analytical solutions and high-accuracy numerical results such as the fourth-order Runge–Kutta and NCFD schemes.

The recursive and symbolic nature of the method not only provides insight into the solution structure but also allows for direct parametric sensitivity analysis. Moreover, the method is computationally efficient and well-suited for implementation in symbolic computation platforms such as Maple, Mathematica, and Python/SymPy.

While this work does not include a formal convergence proof, the empirical results across a wide spectrum of problems consistently demonstrate convergence and accuracy. A full theoretical analysis of convergence properties is a subject for future work.

## Future Work

Future extensions of the method may include:

- Studying the possibility of applying RILT to stochastic differential equations.
- Studying the possibility of extending RILT to handle real powers of both $x$ and $\ln(x)$.
- Application of RILT to IVPs with irregular singularities.
- Application of RILT to IVPs with complex order derivatives.
- Application of RILT to boundary value problems and partial differential equations.
- Application of RILT to eigenvalue problems.

- Symbolic sensitivity analysis and perturbation schemes based on the generated series.
- Acceleration techniques for global convergence using rational approximants or orthogonal basis projections.

**AI Disclosure**

The author used AI tools (specifically OpenAI's ChatGPT) solely for assisting with literature scanning and linguistic polishing. No AI assistance was used in formulating the theoretical framework, derivations, algorithms, or numerical implementations.